\documentclass[12pt]{article}
\ifx\pdftexversion\undefined
  \usepackage[dvips]{graphicx}
\else
  \usepackage[pdftex]{graphicx}
\fi

\usepackage{amsmath}
\usepackage{amssymb}
\usepackage{amsfonts}
\usepackage{amsthm}
\usepackage{epsfig}

\theoremstyle{definition}

\newtheorem{theorem}{Theorem}

\begin{document}

\makebox[1cm]{} 
\vspace{4cm}

\begin{center}
{\Large On a free boundary problem for an } \\

{\Large American put option under the CEV process} 

\end{center}

\begin{center}
{Miao Xu and Charles Knessl}
\end{center}

\begin{center}
{ \small Department of Mathematics, Statistics and Computer Science} \\
{ \small University of Illinois at Chicago} \\
{ \small 851 South Morgan Street} \\
{ \small Chicago, IL 60607-7045}
\end{center}

\begin{center}
{\small e-mail: mxu6@uic.edu, knessl@uic.edu}
\end{center}

\begin{abstract}

We consider an American put option under the CEV process. This corresponds to a free boundary problem for a PDE. We show that this free boundary satisfies a nonlinear integral equation, and analyze it in the limit of small $\rho$ = $2r/ \sigma^2$, where $r$ is the interest rate and $\sigma$ is the volatility. We use perturbation methods to find that the free boundary behaves differently for five ranges of time to expiry.

\end{abstract}


\newpage
\pagestyle{plain}
\setcounter{page}{1}

\section{Introduction}

The pricing and hedging of options has its origins in the Nobel prize winning work of Black, Scholes, and Merton $\cite{cm}$, who assume that the price of an underlying asset $S(t)$ follows a geometric Brownian motion with constant volatility. The price $C(S,t)$ of a European call option at time $t$ for an asset with price $S$, strike $K$, and expiry $T$ is then readily established, and is presented in terms of the normal distribution function. However, there is sufficient empirical evidence $\cite{dl}$ to suggest that in many cases the assumption of constant volatility does not match well to the observed market data. Rather, evidence points out that the implied volatility, which is obtained by equating the model price of an option to its market price and solving for the unknown volatility parameter, varies with the strike price across a wide range of markets. This phenomenon is known as the volatility smile or frown, depending on the shape of the curve, and is not captured by the Blacks-Scholes model with a constant volatility. As a result, there have been various ideas as to how to modify and extend the basic Black-Scholes framework, to account for this phenomenon. One of these is the constant elasticity of variance (CEV) diffusion model, which was introduced by Cox and Ross $\cite{cr}$ in the context of European options. Unlike Black-Scholes, the CEV model is capable of reproducing the volatility smile.  

Other work on European options under a CEV process include Davydov and Linetsky $\cite{dl}$, Hu and Knessl $\cite{hk}$ and Lo, et. al. $\cite{lyh}$. However, there exists little or no analytic work for the valuation of American options under a CEV process. The analysis of these options are more difficult than the corresponding European options in that the American options may be exercised prior to the expiration dates. Mathematically the American options lead to partial differential equations (PDE) with free boundaries, which can only rarely be solved exactly.  In this paper, we apply asymptotic analysis to a CEV model to examine the behavior of the free boundary under different scaling regimes for the time to expiry, in the limit of small $\rho$ = $2r / \sigma^2$, where $r$ is the interest rate, and $\sigma$ is the volatility. This limit has a small interest rate and/or large volatility, and is of particular relevance to the financial status of the current economy. We will employ singular perturbation methods, including matched asymptotic expansions. The main result is the derivation of a nonlinear integral equation that is satisfied by the free boundary, from which we shall analyze its asymptotic structure for five different ranges of time. The main results are summarized in section 2 and derived in section section 3.

Asymptotic analysis and singular perturbation methods have been recently employed in the context of both European and American options,  and this work includes Knessl $\cite{ck1, ck2}$, Howison $\cite{sh}$, Kuske and Keller $\cite{kk}$, Addison, et al. $\cite{ahk}$, Evans, et al. $\cite{ekk}$, Fouque, et al. $\cite{fps}$, and Widdicks, et al. $\cite{wd}$.


\section{Problem Statement and Summary of Results}

We let $P(S,T_0)$ denote the price of an American put option for an asset with price $S$ at some time $T_0$ prior to expiry $T_F$. We assume that $S$ satisfies the stochastic differential equation 
\begin{equation}
dS = \mu S \ dt + \sigma \sqrt{S} \ dW_t. 
\end{equation}
where $W_t$ is a standard Brownian motion, $\sigma$ is the volatility of the underlying asset, and $\mu = r$ is the risk-free interest rate. We note that unlike Black-Scholes, this model only guarantees non-negativity of $S$ ($S \geq 0$), so the chance of absorption at 0, i.e., bankruptcy, occurs with positive probability.  

Introducing the new variables 
\begin{equation}
t = \frac{\sigma^2}{2}(T_F - T_0), \ \  \rho = \frac{2r}{\sigma^2}, 
\end{equation}
we find that $P$ satisfies the following boundary value problem
\begin{equation}
\label{pde1}
P_t = S P_{SS} + \rho SP_S - \rho P; \ t > 0 , \ S > \alpha(t) 
\end{equation}
\begin{equation}
 P(S,0) = \max(K-S,0)
\end{equation}
\begin{equation}  
 P(\alpha(t), t) = K - \alpha(t)
\end{equation}
\begin{equation}
 P_S(\alpha(t), t) = -1
\end{equation}
\begin{equation} 
\label{bc1} 
 P(0,t) = 0
\end{equation}
where $\alpha(t)$ is the free boundary in the new time variable.  We also have $P(S,t) = K - S$ for $0 < S < \alpha(t)$, and $\alpha(0) = K$. For $S \leq \alpha(t)$ the option should be exercised, and for $S > \alpha(t)$ it should be held. 

We convert $(\ref{pde1})- (\ref{bc1})$ into an integral equation by first making a change in coordinates, letting 
\begin{equation}
P(S,t) = K - S + \tilde{P}(V,t),  \ V = S - \alpha(t)   
\end{equation}
where $V$ $\geq$ 0. Then $\tilde{P}$ satisfies the PDE
\begin{equation}
\label{pde2}
\tilde{P}_t - \alpha'(t) \tilde{P}_V = [V + \alpha(t)]\tilde{P}_{VV} + \rho [V + \alpha(t)]\tilde{P}_V - \rho K - \rho \tilde{P}; \ V, t > 0
\end{equation}
with the initial and boundary conditions
\begin{equation}
\tilde{P}(V,0) = V 
\end{equation}
\begin{equation}
\label{bc2}
\tilde{P}(0,t) = \tilde{P}_V (0,t) = 0.
\end{equation}
We introduce the Laplace transform 
\begin{equation}
\label{lp1}
Q(\theta, t) = \int_{0}^{\infty} e^{-\theta V} \tilde{P}(V,t) \ dV.
\end{equation}
Using $(\ref{lp1})$  in $(\ref{pde2})$  and $(\ref{bc2})$ then yields 
\begin{equation}
\label{pde3}
Q_t + (\theta^2 + \rho \theta) Q_\theta = [\alpha'(t) \theta + (\theta^2 + \rho \theta) \alpha(t) - (2\theta + 2\rho)] Q - \frac{\rho K}{\theta}.
\end{equation}
with the initial condition
\begin{equation}
\label{bc3}
Q(\theta, 0) = \frac{1}{\theta^2}.
\end{equation}
Using the method of characteristics, it can be shown that the only acceptable solution to $(\ref{pde3})$ is 
\begin{equation}
\label{pde4}
Q(\theta,t)  =  \frac{K \rho}{\theta^2} e^{\alpha(t)\theta} \int_{\theta / \rho}^{\infty} \frac{1}{z+1}  \exp \left[-\rho z \alpha \left(t+\rho^{-1} \log \left(\frac{\theta+\rho}{\theta}\frac{z}{z+1} \right) \right) \right]  \ dz.
\end{equation}
The next result readily follows. 
\begin{theorem} 

The option price $P(S,t)$ for the CEV model has the integral representation

\begin{equation}
P(S,t) = K - S  + \frac{1}{2 \pi i} \int_{Br} {e^{\theta  V} Q(\theta, t)} \ d\theta,
\end{equation}
where $\Re(\theta) > 0 $ on the Bromwich contour, and $Q(\theta, t)$ is given by $(\ref{pde4})$.     

Moreover, after setting $t = 0$ and using $\alpha(0)=K$ and $(\ref{bc3})$ in $(\ref{pde4})$, it follows that $\alpha(t; \rho)$ satisfies the nonlinear integral equation (IE): 

\begin{equation}
\label{IE}
\frac{e^{-K \theta }}{K \rho} = \int_{\theta / \rho}^{\infty} \frac{1}{z+1} \exp \left[-\rho z \alpha \left(\rho^{-1} \log \left(\frac{\theta+\rho}{\theta}\frac{z}{z+1} \right) \right) \right] \ dz. 
\end{equation}

\end{theorem}

\vspace{0.3in}

In the next section we use asymptotic methods to analyze this IE for five different scales of time $t$, in the limit of small $\rho$. We let $\rho = e^{-\lambda}$ so that $\lambda = - \log \rho \rightarrow \infty$. The final results for the free boundary $\alpha(t; \rho)$ are listed below, and we sketch the derivations in section 3. \\

\begin{enumerate}
 \renewcommand{\labelenumi}{(\roman{enumi})}
   \item $t = \omega / \lambda = O(\lambda^{-1}), \  0 < \omega < K$: 
     \begin{equation}
     \label{eq1}
 \alpha(t; \rho) = (\sqrt{\omega} -\sqrt{K})^2 + \frac{\log \lambda}{\lambda}\frac{\omega - \sqrt{K \omega}}{2} +   \frac{1}{\lambda}    \frac{\sqrt{K \omega} -    \omega}{2} \log \left(\frac{4 \pi K^2 \omega}{K - \       \sqrt{K \omega}}   \right) +      o(\lambda^{-1}) 
     \end{equation} 
   \item $t = K / \lambda + O(\lambda^{-2}), \ \lambda^2 t - \lambda K = \lambda( \omega - K) = \Lambda$:
     \begin{equation}
     \label{eq2}     
   \alpha(t; \rho) = \frac{1}{\lambda^2} \mathcal{F}(\Lambda), 
     \end{equation}
where $\mathcal{F}(\cdot)$ satisfies the nonlinear IE
     \begin{equation}
     \label{eq3}
      \frac{e^{-K\nu}}{K} = \int_{0}^{\infty} \frac{1}{\xi} \exp \left[-\xi F \left(-\nu K^2 - \frac{1}{\xi} \right) \right] \ d\xi.
     \end{equation}  
For $\Lambda \rightarrow \pm \infty$, we have 
     \begin{equation}
     \label{eq4}
      \mathcal{F}(\Lambda) \sim \frac{\Lambda^2}{4K} + \frac{\Lambda}{4} \log( -\Lambda) - \frac{\Lambda}{4} \log(8 \pi K^3),\   \Lambda \rightarrow -\infty 
       \end{equation}
      \begin{equation}
      \label{eq5}
     \mathcal{F}(\Lambda) \sim \Lambda e^{-\gamma} \exp \left[-\frac{1}{K}\exp \left(\frac{\Lambda}{K} \right) \right], \ \Lambda \rightarrow + \infty 
     \end{equation}  
 where $\gamma$ is the Euler constant.    
   \item $t = \omega / \lambda = O(\lambda^{-1}), \ K < \omega < \infty$:
     \begin{equation}
     \label{eq6}
      \alpha(t; \rho) \sim  \frac{\omega - K}{\lambda} e^{-\gamma} \exp \left(-\frac{1}{K} \rho^{K / \omega - 1} \right);
      \end{equation}
   \item  $t = O(1), \ 0 < t < \infty$: 
     \begin{equation}
     \label{eq7}
     \alpha(t; \rho) \sim te^{-\gamma}\exp \left[-\left(\frac{1}{2} + \frac{1}{\rho K} \right) e^{-K/t} \right];
     \end{equation}  
    \item  $t = v / \rho = O(\rho^{-1}), \ v > 0$:
     \begin{equation}
     \label{eq8}
      \alpha(t; \rho) \sim \frac{1}{\rho} e^{-\gamma} \exp \left(\frac{1}{e^{v}-1} \right)(1-e^{-v}) \exp \left(-\frac{1}{\rho K} \right).
     \end{equation}
\end{enumerate}


We note that in four of the five cases the expression for $\alpha(t; \rho)$ is completely explicit, and only in case (ii) must we solve a nonlinear IE, which is somewhat simpler than the one in $(\ref{IE})$. We can easily compute $P(S,t)$ as $t \rightarrow \infty$, which corresponds to the perpetual American option, where the problem reduces to solving an ordinary differential equation. Setting $P(S,\infty) = P^\infty(S)$ and using $\alpha(\infty)$ to denote the limiting value of the free boundary, we obtain from $(\ref{pde1})-(\ref{bc1})$   

\begin{equation}
P^\infty(S) =  Ke^{\rho \alpha(\infty)} \int_{1}^{\infty} \frac{1}{z^2}e^{-z \rho S} \ dz,  
\end{equation}
where $\alpha(\infty)$ satisfies
\begin{equation}
\label{aper}
K \rho \int_{1}^{\infty} \frac{1}{z}e^{-\rho \alpha(\infty) z} \ dz = e^{-\rho \alpha(\infty)}.
\end{equation}
For $\rho \rightarrow 0$ we have 
\begin{equation}
\label{fbper}
\alpha(\infty) = \frac{1}{\rho} e^{-\gamma} \exp \left(-\frac{1}{\rho K} \right) [1 + O(\rho)],  
\end{equation}
which is exponentially small.


\section{Analysis}

\subsection{Analysis for $t = \omega / \lambda$, $0 < \omega < K $}

We first examine $(\ref{IE})$ on the $t = O(\lambda^{-1})$ scale, for small $\rho$. Recalling that $\lambda = -\log \rho$, we let 
\begin{equation}
\theta = \lambda \beta, \ z = \frac{\lambda(\beta + x)}{\rho}, \ \alpha(t; \rho) \sim \alpha_0(\omega; \rho) 
\end{equation}
for $\omega = (-\log \rho)t = O(1)$. Then $(\ref{IE})$ can be approximated by 
\begin{equation}
\label{ie1}
\frac{e^{\lambda}}{K} = \int_{0}^{\infty} \frac{1}{\beta + x} e^{\lambda \Phi(x; \beta, \rho)}  [1 + O(e^{-\lambda})] \ dx,
\end{equation}    
where $\Phi(x; \beta, \rho) = K \beta - (\beta + x) \alpha_0(\frac{x}{\beta(x + \beta)}; \rho)$. For large $\lambda$ and fixed $\beta$, we evaluate the right hand side of $(\ref{ie1})$ by an implicit form of the Laplace method, assuming for now that there is a saddle point where
\begin{equation}
\label{saddle}
\frac{\partial \Phi}{\partial x} = -\alpha_0 \left(\frac{x}{\beta(x + \beta)} \right) - \frac{1}{x + \beta}\alpha_0' \left(\frac{x}{\beta(x + \beta)} \right) = 0.
\end{equation}
Let us denote $x = x_*(\beta)$ as the solution to $(\ref{saddle})$. It follows that at $x = x_*$, $\Phi \sim 1$ so that 
\begin{equation}
\label{spe}
1 = K\beta - (\beta + x_*)\alpha_0 \left(\frac{x_*}{\beta(x_* + \beta)} \right).
\end{equation} 
Now let $\omega = \frac{x_*}{\beta(x_* + \beta)}$. Then from $(\ref{saddle})$ we have $\beta = \frac{\alpha_0'(\omega)}{\omega \alpha_0'(\omega) - \alpha_0(\omega)}$ which we use to eliminate $\beta$ in $(\ref{spe})$ to obtain the ODE

\begin{equation}
[1-\alpha_0'(\omega)][\omega \alpha_0'(\omega) - \alpha_0(\omega)] = K \alpha_0'(\omega)
\end{equation}
Rewriting this as
\begin{equation} 
\omega \alpha_0'(\omega) - \alpha_0(\omega) = \frac{K \alpha_0'(\omega)}{1-\alpha_0'(\omega)},
\end{equation}
we recognize this as the Clairaut equation. The solutions consist of a one-parameter family of lines and the singular solution 
\begin{equation}
\label{wsol}
\alpha_0(\omega) = (\sqrt{\omega} - \sqrt{K})^2
\end{equation}
which is the envelope of this family. The linear solutions $\alpha_0(\omega) = \omega C - \frac{KC}{1-C}$ must be rejected, since these lead to $\alpha_0(0) \neq K$. The above analysis applies only for $0 < \omega < K$, since the solution $(\ref{wsol})$ vanishes as $\omega$ approaches $K$. Hence we expect different asymptotics for $\omega \approx K$. 

We next analyze some higher order terms in the expansion of $\alpha_0$. We evaluate $(\ref{ie1})$ by using the Laplace method, which gives  

\begin{equation}
\label{laplace}
\frac{e^{\lambda}}{K} = \frac{1}{\beta + x_*} \sqrt{\frac{2\pi}{-\lambda \Phi_{xx}(x_*; \beta, \rho)}} e^{\lambda \Phi(x_*; \beta, \rho)}[1 + O(\lambda^{-1})],
\end{equation}
and expand $\alpha_0$ as 
\begin{equation}
\label{exp1}
\alpha_0(\omega; \rho) = \alpha_0(\omega) + \frac{\log \lambda}{\lambda} \alpha_1(\omega) + \frac{1}{\lambda} \alpha_2(\omega) + o(\lambda^{-1}).
\end{equation}  
In order to balance the two sides of $(\ref{laplace})$, we need $\alpha_1$ to cancel the $\sqrt{1 / \lambda}$ factor. Hence,
\begin{equation}
\label{la1}
\frac{1}{\sqrt{\lambda}}\exp \left[-(\beta + x_*)(\log \lambda) \ \alpha_1 \left(\frac{x_*}{\beta(x_* + \beta)} \right) \right] = 1.
\end{equation} 
Writing $(\ref{la1})$ in terms of $\omega$ we obtain 
\begin{equation}
\label{wsol2}
\alpha_1(\omega) = \frac{1}{2} (\omega - \sqrt{K\omega}).
\end{equation} 
To find the third order term $\alpha_2$, we balance the $O(1)$ terms in $(\ref{laplace})$, so that
\begin{equation}
\frac{1}{K} = \frac{1}{\beta + x_*} \sqrt{\frac{{2\pi}}{-\Phi_{xx}}} \exp \left[-(\beta + x_*) \alpha_2 \left(\frac{x_*}{\beta(x_* + \beta)} \right)\right].
\end{equation}
It can be shown that $\Phi_{xx}(x_*; \beta, \rho) \sim -\frac{1}{2} K^{\frac{1}{2}} \beta^{\frac{3}{2}}x_*^{-\frac{3}{2}}(x_* + \beta)^{-\frac{3}{2}}$ and then
\begin{equation}
\label{wsol3}
\alpha_2(\omega) = \frac{\sqrt{K\omega} - \omega}{2} \log \left(\frac{4\pi K^2 \omega}{K - \sqrt{K\omega}} \right).
\end{equation}
With $(\ref{wsol})$, $(\ref{exp1})$, $(\ref{wsol2})$, and $(\ref{wsol3})$ we have established $(\ref{eq1})$.

\subsection{Analysis for $t = \omega / \lambda$, $\omega \approx K $}

We return to $(\ref{IE})$ and introduce the scaling
\begin{equation}
\label{scale1}
\theta = \lambda \beta, \ \beta = \frac{1}{K} + \frac{\nu}{\lambda}, \ \omega = \lambda t = K + \frac{\Lambda}{\lambda}
\end{equation}
with
\begin{equation}
\label{scale2}
\alpha(t; \rho) = \frac{1}{\lambda^2} \mathcal{F}(\Lambda; \rho) = \frac{1}{\lambda^2} \mathcal{F} \left(\lambda^2 \left(t-\frac{K}{\lambda} \right); \rho \right).
\end{equation}
Then we have $e^{-\theta K} = \rho e^{-K \nu}$. Also by setting $z = (\theta + y) / \rho $ in $(\ref{IE})$ this equation becomes 
\begin{equation}
\label{ie2}
\frac{1}{K} e^{-K \nu} = \int_{0}^{\infty} \frac{1}{\theta + y + \rho} \exp \left[-(\theta + y) \alpha \left(\frac{1}{\theta} - \frac{1}{\theta + y} + O(\rho); \rho \right) \right] \ dy.
\end{equation}
With the scaling in $(\ref{scale1})$ and $(\ref{scale2})$ and the fact that $\rho = e^{-\lambda}$ is exponentially small, we obtain 
\begin{align}
\lambda^2 \alpha \left(\frac{1}{\theta} - \frac{1}{\theta + y} + O(e^{-\lambda}); \rho \right) & = \lambda^2 \alpha \left(\frac{K}{\lambda} - \frac{1}{\lambda^2} \left(\nu K^2 + \frac{1}{\xi} \right) + o(\lambda^{-2}); \rho \right)\notag \\
& \sim \mathcal{F} \left(-\nu K^2 - \frac{1}{\xi} \right)
\end{align}
where $\mathcal{F}(\Lambda)$ is the leading term in an expansion of $\mathcal{F}(\Lambda; \rho)$. Then scaling $y = \lambda^2 \xi$ in $(\ref{ie2})$ and letting $\rho \rightarrow 0$ $(\lambda \rightarrow \infty)$ we obtain the limiting IE in $(\ref{eq3})$. It does not seem possible to solve $(\ref{eq3})$ explicitly for $\mathcal{F}(\Lambda)$. But we can infer the behavior as $\Lambda \rightarrow -\infty$ by evaluating the integral in $(\ref{eq3})$ by an implicit Laplace type expansion, similarly to what we did in section 3.1. This will verify the asymptotic matching between the $\omega$-scale (for $\omega < K$) and the $\Lambda$-scale, and lead to $(\ref{eq4})$. Now consider the limit $\Lambda \rightarrow +\infty$. For $|\nu| < 0$, we rewrite $(\ref{eq3})$  as 

\begin{equation}
\label{inte}
\begin{split}
\frac{e^{K |\nu|}}{K}  & = \int_{\lvert \nu \rvert K^2}^{\infty} \frac{1}{\eta} \exp \left[-\frac{1}{\eta} \mathcal{F}(\lvert \nu \rvert K^2 - \eta) \right] \ d\eta  \\ & + \int_{0}^{1} \frac{1}{u} \left \{ \exp \left[-\frac{1}{u} \frac{\mathcal{F}(\lvert \nu \rvert K^2 (1-u))}{\lvert \nu \rvert K^2} \right] - \exp \left[-\frac{1}{u} \frac{\mathcal{F}(\lvert \nu \rvert K^2)}{\lvert \nu \rvert K^2} \right] \right \}  \\ & + 
\int_{\frac{\mathcal{F}(\lvert \nu \rvert K^2)}{\lvert \nu \rvert K^2}}^{\infty} \frac{e^{-v}}{v} \ dv. \hspace{1in}  
\end{split}
\end{equation}  										
Here we broke up the integral over $(0, \infty)$ into the two ranges $(0, \lvert \nu \rvert K^2)$ and $(\lvert \nu \rvert K^2, \infty)$ and made some elementary substitutions. Now, for $\Lambda \rightarrow  - \infty$ we have $\mathcal{F}(\Lambda) \sim \frac{\Lambda^2}{4K}$ so that the first integral in the right hand side of $(\ref{inte})$ will vanish as $\nu \rightarrow -\infty$. If $\mathcal{F}(\Lambda) \rightarrow 0$ as $\Lambda \rightarrow + \infty$ the second integral in $(\ref{inte})$ will also vanish, and the third may be approximated by using
\begin{equation}
\label{expint}
\int_{\varepsilon}^{\infty} \frac{e^{-v}}{v} \ dv =  - \log \varepsilon - \gamma + O(\varepsilon),\;\; \varepsilon \rightarrow 0^+. 
\end{equation} 
Hence $(\ref{inte})$ can be replaced by the asymptotic relation 
\begin{equation}
\frac{e^{K \lvert \nu \rvert}}{K} \sim -\log \left[\frac{\mathcal{F}(\lvert \nu \rvert K^2)}{\lvert \nu \rvert K^2} \right] - \gamma
\end{equation}
which upon exponentiation leads to the asymptotic result given in $(\ref{eq5})$, for $\mathcal{F}(\Lambda)$ as $\Lambda \rightarrow \infty$. 

\subsection{Analysis for $t = \omega / \lambda$, $K < \omega < \infty $}
In the remaining time ranges, $\alpha(t; \rho)$ will be exponentially small as $\rho = e^{-\lambda} \rightarrow 0$, and our analysis of $(\ref{IE})$ will rely heavily on the asymptotic form in $(\ref{expint})$. 
We let $z = Z / \rho$ in $(\ref{IE})$ to obtain 
\begin{equation}
\label{ie3}
\frac{e^{-K \theta}}{K \rho} = \int_{\theta}^{\infty} \frac{1}{Z+\rho} \exp \left[-Z \alpha \left(\rho^{-1} \log \left(\frac{\theta + \rho}{\theta} \frac{Z}{Z + \rho} \right) \right) \right] \ dZ.
\end{equation}
Now we scale $Z = \lambda z_*$ and $\theta = \lambda \theta_*$, let $\alpha(t; \rho) = \tilde{\alpha}(\lambda t; \rho)$ and note that in sections 3.1 and 3.2 we have already characterized $\tilde{\alpha}(\lambda t; \rho)$ for $\lambda t = \omega < K$ and $\omega \sim K$. We also simplify the argument of $\alpha(\cdot)$ in $(\ref{ie3})$ using 
\begin{equation}
\alpha \left(\frac{1}{\rho} \left[\log \left(1 + \frac{\rho}{\theta} \right) - \log \left(1 + \frac{\rho}{Z} \right) \right] \right) = \alpha \left(\frac{1}{\theta} - \frac{1}{Z} + O(\rho) \right) = \tilde{\alpha} \left(\frac{1}{\theta_*} - \frac{1}{Z_*} + O(e^{-\lambda} \lambda) \right).
\end{equation}
When $Z = \theta$ we have $z_* = \theta_*$ and we rewrite the integral in $(\ref{ie3})$ by splitting the range of integration into $z_* \in (\theta_*, \theta_* / [1 - K \theta_*])$ and $z_* \in (\theta_* / [1 -K \theta_*], \infty)$, thus obtaining 
\begin{equation}
\frac{e^{-K \lambda \theta_*}}{K \rho} \sim \left( \int_{\theta_*}^{\frac{\theta_*}{1-K \theta_*}} + \int_{\frac{\theta_*}{1-K \theta_*}}^{\infty} \right) \left( \frac{1}{z_*} \exp \left[-\lambda z_* \tilde{\alpha} \left( \frac{1}{\theta_*} - \frac{1}{z_*} \right) \right] \right) \ dz_*.
\end{equation}
In the first range $\tilde{\alpha}(\omega) \sim (\sqrt{K} - \sqrt{\omega})^2$ and the first integral will be $o(1)$ as $\lambda \rightarrow \infty$, since $\theta_*^{-1} - z_*^{-1} \leq K$ when $z_* \leq \theta_* / [1 - K \theta_*]$. In the second integral $\tilde{\alpha}$ will be exponentially small and the main contribution will come from very large values of $z_*$, where roughly $z_* = O(\tilde{\alpha}^{-1})$. Then we write $\tilde{\alpha}(\theta_*^{-1} - z_*^{-1}) \sim \tilde{\alpha}(\theta_*^{-1})$ and using $(\ref{expint})$ we conclude that
\begin{equation}
\label{ie4}
\frac{e^{-K \lambda \theta_*}}{K \rho} \sim \int_{\frac{\theta_*}{1 - K \theta_*}}^{\infty} \frac{1}{z_*} \exp \left[-\lambda z_* \tilde{\alpha} \left(\frac{1}{\theta_*} \right) \right] \ dz_* \sim - \log \left[\lambda \tilde{\alpha} \left(\frac{1}{\theta_*} \right) \right] - \gamma  - \log \left[\frac{\theta_*}{1 - K \theta_*} \right], 
\end{equation}
with an error that is $o(1)$ as $\lambda \rightarrow \infty$. Then exponentiating $(\ref{ie4})$ and replacing $\theta_*$ by $\omega^{-1}$ we obtain the asymptotic result in $(\ref{eq6})$.

For $\omega \rightarrow K$ we note that $\rho^{\omega / K - 1} = \rho^{-1} e^{-K / t} = \rho^{-1} \exp \left[-\frac{\lambda K}{K + \frac{\Lambda}{\lambda}} \right]  \\
= \rho^{-1} \exp \left[-\lambda + \frac{\Lambda}{K} + O(\lambda^{-1}) \right] \sim \exp \left(\frac{\Lambda}{K} \right)$ 
and $(\omega - K) / \lambda = \Lambda / \lambda^2$, which can be used to verify the asymptotic matching between the $\Lambda$-scale and the $\omega$-scale for $\omega > K$, in the intermediate limit where $\omega \downarrow K$ and $\Lambda \rightarrow \infty$.

\subsection{Analysis for $t = O(1), 0 < t < \infty $}
Next we consider times $t = O(1)$. We scale $z = \theta w / \rho$. Since we again expect $\alpha(t; \rho)$ to be very small we assume a ``WKB-type'' ansatz of the form  
\begin{equation}
\label{wkb}
\alpha(t; \rho) \sim g(t) \exp \left[-\frac{1}{\rho} f(t) \right].
\end{equation}
Expanding $\alpha(t; \rho)$ in $(\ref{IE})$  
for fixed $\theta$ and $\rho \rightarrow 0$, and noting that 
\begin{equation*}
\frac{1}{\rho} \left[\log \left(1 + \frac{\rho}{\theta} \right) - \log \left(1 + \frac{\rho}{\theta w} \right) \right] = \frac{1}{\theta} - \frac{1}{\theta w} - \frac{1}{2} \frac{\rho}{\theta^2} + O \left(\rho^2, \frac{\rho}{w^2} \right),
\end{equation*}
we have
\begin{align*}
-\theta w \alpha \left( \frac{1}{\rho} \left[\log \left(1 + \frac{\rho}{\theta} \right) - \log \left(1 + \frac{\rho}{\theta w} \right) \right] \right) & \sim -\theta w g \left(\frac{1}{\theta} \right) \exp \left[-\frac{1}{\rho} f \left(\frac{1}{\theta} - \frac{1}{\theta w} - \frac{1}{2} \frac{\rho}{\theta^2} \right) \right] \\
& = -\theta w g \left(\frac{1}{\theta} \right) \exp \left[-\frac{1}{\rho} f \left(\frac{1}{\theta} \right) + \frac{1}{2 \theta^2} f' \left(\frac{1}{\theta} \right) + O(\rho) \right]
\end{align*}
Here we also used $f(\theta^{-1} - (\theta w)^{-1}) \sim f(\theta^{-1})$, since $w$ will be scaled to be exponentially large. Then setting $\varepsilon = \theta g \left(\frac{1}{\theta} \right) \exp \left[-\frac{1}{\rho} f \left(\frac{1}{\theta} \right) \right] \exp \left[\frac{1}{2 \theta^2} f' \left(\frac{1}{\theta} \right) \right]$, 
scaling $w = \varepsilon^{-1} u$ and using $(\ref{expint})$, $(\ref{IE})$ asymptotically becomes
\begin{equation}
\label{asy1}
\frac{e^{-K \theta }}{K \rho} \sim \frac{1}{\rho} f \left(\frac{1}{\theta}\right) -\gamma - \log \left[\theta g \left(\frac{1}{\theta} \right) \right] - \frac{1}{2\theta^2}f' \left(\frac{1}{\theta} \right) + o(1)
\end{equation} 
From the $O(\rho^{-1})$ terms in $(\ref{asy1})$ we conclude that $f \left(1/ \theta \right) = K^{-1} e^{-K \theta}$ and then the $O(1)$ terms determine $g(\cdot)$ from 
\begin{equation*}
\label{match}
\theta g \left(\frac{1}{\theta} \right) = e^{-\gamma} \exp \left[-\frac{1}{2 \theta^2} f' \left(\frac{1}{\theta} \right) \right]. 
\end{equation*}
The above along with $(\ref{wkb})$ establishes the asymptotic result in $(\ref{eq7})$. 
The asymptotic matching between $(\ref{eq6})$ and $(\ref{eq7})$ is immediate, since $\rho^{K/ \omega - 1} = \rho^{-1} e^{-K / t}$, and $(\omega - K) / \lambda \sim \omega / \lambda = t$ as $\omega \rightarrow \infty$.

\subsection{Analysis for $t = v / \rho  = O(\rho^{-1}), v > 0$}
We assume that time to expiry for the option is large, with $t = v / \rho = O(\rho^{-1})$. On this time scale we assume that 
\begin{equation}
\label{vwkb}
\alpha(t; \rho) \sim \frac{1}{\rho} \exp \left(-\frac{1}{\rho K} \right) A(v),
\end{equation}
where $A(\cdot)$ will be determined from $(\ref{IE})$.
After scaling $\theta = \rho W$, $(\ref{IE})$ becomes 
\begin{equation}
\label{vie}
\frac{e^{-K \rho W}}{K \rho} \sim \int_{W}^{\infty} \frac{1}{z+1} \exp \left[-ze^{-\frac{1}{\rho K}}A \left(\log \left[\frac{W+1}{W} \frac{z}{z+1} \right] \right) \right] \ dz.
\end{equation}   
The major contribution to the integral in $(\ref{vie})$ will once more come from large values of $z$, so we approximate 
\begin{equation}
A \left(\log \left[\frac{W + 1}{W} \frac{z}{z + 1} \right] \right) \sim A \left(\log \left[\frac{W + 1}{W} \right] \right),
\end{equation}
and then applying this to $(\ref{IE})$ along with $(\ref{expint})$ leads to 
\begin{align}
\frac{e^{-K \rho W}}{K \rho} & \sim - \log \left\{e^{-\frac{1}{\rho K}} A \left(\log \left[\frac{W + 1}{W} \right] \right) \right\} - \gamma  - \log(W + 1) \\
&= \frac{1}{K \rho} - \gamma - \log \left\{(W + 1) A \left(\log \left[\frac{W + 1}{W} \right] \right) \right\} + o(1).
\end{align}
Then after expanding $e^{-K \rho W} = 1 - K \rho W + O(\rho^2)$ we conclude that 
\begin{equation}
A \left(\log \left[\frac{W + 1}{W} \right] \right) = e^{-\gamma} \frac{1}{W + 1} e^W
\end{equation}  
which determines the function $A(\cdot)$ and establishes $(\ref{eq8})$. 

Finally we verify the asymptotic matching between $(\ref{eq7})$ and $(\ref{eq8})$. For $v \rightarrow 0$ we have $(e^v - 1)^{-1} = v^{-1} - \frac{1}{2} + O(v)$ and $1 - e^{-v} \sim v = \rho t$. For $t \rightarrow \infty$ we have $e^{-K / t} = 1 - \frac{K}{t} + O(t^{-2}) = 1 - \frac{K \rho}{v} + O(\rho^2)$ so that $-(\frac{1}{2} + \frac{1}{\rho K}) e^{-K / t} \sim -v^{-1} + \frac{1}{2}$ and the matching follows. As $v \rightarrow \infty$ we have $A(v) \rightarrow e^{-\gamma}$ and thus the expansion in $(\ref{vwkb})$ agrees with the small $\rho$ expansion of $\alpha(\infty; \rho)$, as given in $(\ref{fbper})$. 


\end{document}